\documentclass[a4paper,12pt]{article}
\usepackage{latexsym,amssymb,amsmath,amsthm}
\usepackage[all]{xy}
\usepackage[T1]{fontenc}
\usepackage{lmodern}

\newtheorem{theorem}{Theorem}
\newtheorem{proposition}[theorem]{Proposition}
\newtheorem{corollary}[theorem]{Corollary}

\newcommand{\N}{{\mathbb N}}
\newcommand{\Z}{{\mathbb Z}}
\newcommand{\Q}{{\mathbb Q}}
\newcommand{\F}{{\mathbb F}}
\def\mapd#1{\Big\downarrow\rlap {$\vcenter{\hbox{$\scriptstyle{{#1}}$}}$}}

\def\ds{\displaystyle}

\newdimen\Itemindent \Itemindent=.9cm

\def\lpfeil{\mathop{\longrightarrow}\limits}
\def\pfeil{\mathop{\rightarrow}\limits}
\def\Hom{\mathrm{Hom}}
\def\lr@iso{\mathrel{\kern6pt\lower2pt\hbox{$\scriptstyle\sim$} \kern-12pt\hbox{$\longrightarrow$}}}
\def\lriso#1{\mathrel{\mathop{\lr@iso}\limits^{#1}}}
\def\cor{{\mathop{\it cor\,}}}
\def\cd{\mathit{cd}}
\def\res{{\mathop{\it res\,}}}
\newcommand{\wingl}{[\![}
\newcommand{\wingr}{]\!]}
\def\Mod{\mathit{Mod}}
\def\tr{\mathit{tr}}

\def\vcd{\mathit{vcd}}
\newcommand{\hatotimes}{\mathop{\hat\otimes}}

\title{\bf\boldmath  Extensions of profinite duality groups}
\author{Alexander Schmidt  and Kay Wingberg}
\date{September 25, 2008}
\begin{document}
\maketitle

\indent

Let $G$ be a profinite group and let  $p$ be a prime number. By  $\Mod_p(G)$ we denote the category of discrete $p$-primary $G$-modules. For $A\in \Mod_p(G)$ and $i\geq 0$, let
$$
D_i(G,A)=\varinjlim_U\, H^i(U,A)^*,
$$
where $\null^*$ is $\Hom(-,\Q_p/\Z_p)$, the direct limit is taken over all open subgroups $U$ of $G$ and the transition maps are the duals of the corestriction maps. $D_i(G,A)$ is a discrete $G$-module in a natural way.  Assume that $n=\cd_p\,G$ is finite. Then the $G$-module
$$
I(G)=\varinjlim_{\nu\in \N}\,D_n(G,\Z/p^\nu\Z)
$$
is called the {\bf\boldmath dualizing module} of $G$ at $p$. Its importance lies in the functorial
isomorphism
$$
H^n(G,A)^*\cong\Hom_G(A,I(G))
$$
for all $A\in\Mod_p(G)$. This isomorphism is induced by the cup-products ($V\subseteq U$)
\[
H^n(G,A)^* \times \null_{p^\nu} A^U \longrightarrow H^n(V,\Z/p^\nu\Z)^*,\ (\phi,a)\longmapsto \big(\alpha \mapsto \phi(\cor_G^V(\alpha \cup a))\big)
\]
by passing to the limit over $\nu$ and $V$, and then over $U$.
The identity-map of $I(G)$ gives rise to the homomorphism
$$
\tr :H^n(G,I(G))\lpfeil \Q_p/\Z_p\,,
$$
called the {\bf trace map}.

\medskip
The profinite group $G$ is called a {\bf\boldmath duality group at $p$ of dimension $n$} if
for all $i\in\Z$ and all finite $G$-modules $A\in\Mod_p(G)$,
the cup-product and the trace map
\smallskip\noindent
$$
H^i(G,\Hom(A,I(G)) \times H^{n-i}(G,A)\lpfeil^\cup
H^n(G,I(G))\lpfeil^{\tr}\Q_p/\Z_p
$$

\smallskip\noindent
yield an isomorphism

\smallskip\noindent
$$
H^i(G,\Hom(A,I(G)))\cong H^{n-i}(G,A)^*.
$$

\bigskip\noindent
{\bf Remark:} In \cite{V}, J.-L.~Verdier  used the name {\bf\boldmath strict Cohen-Macaulay at $p$} for what we call  a profinite duality group at $p$ here.  In \cite{48}, A.~Pletch defined   $D_p^n$-groups (and called them duality groups at $p$ of dimension~$n$). The $D_p^n$-groups of Pletch are exactly the duality groups at $p$ (in our sense) which, in addition,  satisfy the following finiteness condition:

\medskip\noindent
$FC(G,p)$: \,\,{\it $H^i(G,A)$ is finite  for all finite $A\in \Mod_p(G)$ and for all $i\geq 0$}.

\medskip\noindent
Since any finite, discrete $G$-module is trivialized by an open subgroup $U$ of $G$,  condition $FC(G,p)$ can also be rephrased in the form:

\medskip\noindent
$FC(G,p)$: \,\,{\it $H^i(U,\Z/p\Z)$ is finite for all open subgroups $U$ of $G$ and
all $i\geq 0$}.

\bigskip
By a duality theorem due to J. Tate, see \cite{T} Thm.\,3 or \cite{V} Prop.\,4.3 or  \cite{NSW} (3.4.6),
a profinite group $G$ of cohomological $p$-dimension  $n$ is a duality group at $p$ if and only if
$$
D_i(G,\Z/p\Z)=0 \quad \hbox{ for $0\leq i<n$}.
$$
As a consequence we see that every open subgroup of a duality group at $p$ is  a duality group at $p$
as well (of the same cohomological dimension), and if an open subgroup of $G$ is a duality group at $p$ and $\cd_p\,G<\infty$, then $G$
is  a duality group at $p$ of the same cohomological dimension (use \cite{NSW} (3.3.5)(ii)). Furthermore, any profinite group of cohomological $p$-dimension~$1$ is a duality group at $p$.

\medskip
We call a profinite group $G$  {\bf\boldmath virtually a duality group at $p$ of (virtual) dimension
$\vcd_p \,G=n$} if
an open subgroup $U$ of $G$ is a duality group at $p$ of dimension $n$.

\bigskip
The objective of this paper is to give a proof of Theorem~\ref{1} below, which states that the class of duality groups is closed under group extensions $1\to H\to G\to G/H \to 1$ if the kernel satisfies $FC(H,p)$. Weaker forms of  Theorem~\ref{1} were first proved by
A.~Pletch (for $D_p^n$-groups, see \cite{48}\footnote{The proof given by Pletch in \cite{48} is only correct for pro-$p$-groups as the author assumes that finitely generated projective modules over the complete group ring
$\Z_p\wingl G\wingr$ are free.}) and by the second author (for Poincar\'{e}  groups, see \cite{71}).

\begin{theorem}\label{1}
Let
\smallskip\noindent
$$
1 \lpfeil H\lpfeil G\lpfeil G/H\lpfeil 1
$$

\smallskip\noindent
be an exact sequence of profinite groups such that condition $FC(H,p)$ is satisfied.
Then the following assertions hold:

\begin{description}
\item{\rm (i)} If\/ $G$ is a duality group at $p$, then $H$ is  a duality group at $p$ and $G/H$ is virtually a duality group at $p$.
\item{\rm (ii)} If\/ $H$ and $G/H$ are duality groups at $p$, then $G$ is a duality group at $p$.
\end{description}
Moreover, in both cases we have:
$$
\cd_p\,G=\cd_p\,H+\vcd_p\,G/H,
$$

\smallskip\noindent
and there is a canonical $G$-isomorphism
\[
I(G)^\vee\cong I(H)^\vee \; {\hatotimes}_{\Z_p}\; I(G/H)^\vee,
\]
where $\null^\vee$ is the Pontryagin dual and ${\hatotimes}_{\Z_p}$ is the tensor-product in the category of compact $\Z_p$-modules.
\end{theorem}

\noindent
{\bf Remark:} The assumption $FC(H,p)$ is necessary, as the following examples show:
\begin{itemize}
\item[1.] Let $G$ be the free pro-$p$-group on two generators $x,y$ and let $H\subset G$ be the normal subgroup generated by $x$. Then $H$ is free of infinite rank, $G/H$ is free of rank one and $1 \pfeil H\pfeil G\pfeil G/H\pfeil 1$ is an exact sequence in which all three groups are duality groups of dimension one.
\item[2.]  Let $D$ be a duality group at $p$ of di\-men\-sion $2$, $F$ a duality group at $p$ of dimension $1$  and $G=F\ast D$ their free product. The kernel of the projection $G\twoheadrightarrow D$ has cohomological $p$-dimension~$1$, hence is a duality group a $p$ of dimension $1$. The group $G$ has cohomological $p$-dimension~$2$ but is is not a duality group at $p$.
\end{itemize}

\medskip\noindent
In the proof of Theorem~\ref{1}, we make use of the following

\begin{proposition} \label{spectral}
Let
\[
1 \lpfeil H\lpfeil G\lpfeil G/H\lpfeil 1
\]
be an exact sequence of profinite groups. Assume that FC($H,p$) holds. Then there is a spectral sequence of homological type
\[
E^2_{ij}= D_i(G/H,\Z/p\Z) \otimes D_j(H,\Z/p\Z) \Longrightarrow D_{i+j}(G,\Z/p\Z).
\]
\end{proposition}

\begin{proof}
Let $g$ run through the open normal subgroups of $G$.
Then $gH/H \cong g/g\cap H$ runs through the open normal subgroups of $G/H$. For a
$G$-module $A\in\Mod_p(G)$, we consider the Hochschild-Serre spectral
sequence
\smallskip\noindent
$$
E(g,g\cap H,A):\ E_2^{ij}(g,g\cap H,A)=H^i(g/g\cap H, H^j(g\cap H,A))\Longrightarrow H^{i+j}(g,A).
$$

\smallskip\noindent
If $g'\subseteq g$ is another open normal subgroup of $G$, then the
corestriction yields a morphism
$$
\cor:E(g',g'\cap H,A)\lpfeil E(g,g\cap H,A)
$$
of spectral sequences.
The map
$$
E_2^{ij}(g',g'\cap H,A) \lpfeil E_2^{ij}(g,g\cap H,A)
$$
is the composite of the maps
\[
H^i\big(g'/g'\cap H, H^j(g'\cap H,A)\big)
\lpfeil^{\cor_{g\cap H}^{g'\cap H}} H^i\big(g'/g'\cap H,H^j(g\cap H,A)\big)
\]
\[
\lpfeil^{\cor_{g/g\cap H}^{g'/g'\cap H}} H^i\big(g/g\cap H,H^j(g\cap H,A)\big)
\]
\smallskip\noindent
and the map between the limit terms is the corestriction
$$
\cor^{g'}_g:H^{i+j}(g',A)\lpfeil H^{i+j}(g,A).
$$
For $2\leq r\leq\infty$ we set
$$
E^2_{ij}=D^r_{ij}(G,H,A):=\varinjlim_g\,E^{ij}_r(g,g\cap H,A)^*.
$$
As taking duals and direct limits are exact operations, the terms $D^r_{ij}(G,H,A)$, $2\leq r \leq \infty$, establish a homological spectral sequence which converges to $D_n(G,A)$.  If $h$ runs through the open subgroups of $H$ which are normal in $G$, then
the cohomology groups $H^j(h,A)$ are $G$-modules in a natural way. If $g$ is open in $G$ with $g\cap H\subseteq h$, then these groups are  $g/g\cap H$-modules. We see that

\smallskip\noindent
$$
D^2_{ij}(G,H,A)=\varinjlim_{\substack{h\subseteq H\\h\trianglelefteq G}}\
\varinjlim_{\substack{g\subseteq G\\ g\cap H\subseteq h}}\,H^i(g/g\cap H,H^j(h,A))^*,
$$

\smallskip\noindent
where for both limits the transition maps are (induced by) $\cor^*$.
In order to conclude the proof of the proposition, it remains to construct isomorphisms
\[
D^2_{ij}(G, H,\Z/p\Z)\cong D_{i} (G/H,\Z/p\Z) \otimes D_{j}(H,\Z/p\Z)
\]
for all $i$ and $j$. To this end note that all occurring abelian groups are $\F_p$-vector spaces, so that $\null^*$ is $\Hom(-,\F_p)$. Further note that for vector spaces $V,W$ over a field $k$ the homomorphism
\[
V^*\otimes W^* \lpfeil (V\otimes W)^*,\ \phi\otimes \psi \longmapsto \big(v\otimes w\mapsto \phi(v)\psi(w)\big)
\]
is an isomorphism provided that $V$ or $W$ is finite-dimensional.  Let $h$ be an  open subgroup of $H$ which is normal in $G$ and let $g'\subseteq g$ be open subgroups of $G$ such that $g$ acts trivially on the finite group $H^{j}(h,\Z/p\Z)$. Then, by \cite{NSW} (1.5.3)(iv), the diagram
\[
\renewcommand{\arraystretch}{1.6}
\begin{array}{ccc}
 H^{i}(g'/g'\cap H,\Z/p\Z)\otimes H^{j}(h,\Z/p\Z) & \lriso{\cup} & H^{i}\big(g'/g'\cap H, H^{j}(h,\Z/p\Z)\big)\\
 \phantom{xxxmN}\mapd{\cor \otimes \mathit{id}}&&\mapd{\cor}\\
 H^{i}(g/g\cap H,\Z/p\Z)\otimes H^{j}(h,\Z/p\Z) &  \lriso{\cup} & H^{i}\big(g/g\cap H, H^{j}(h,\Z/p\Z)\big)
\end{array}
\renewcommand{\arraystretch}{1}
\]
commutes. For fixed $h$, we therefore obtain isomorphisms

\bigskip
$
\ds D_{i}(G/H,\Z/p\Z) \otimes H^{j}(h,\Z/p\Z)^*
$
\[
\begin{array}{lcl}
\phantom{uuuuuuuuuuuuuuuuuuuu}&\cong&
\ds\big(\varinjlim_g H^{i}(g/g\cap H,\Z/p\Z)^*\big) \otimes H^{j}(h,\Z/p\Z)^* \\
&\cong&\ds\varinjlim_g \; H^{i}(g/g\cap H,\Z/p\Z)^* \otimes H^{j}(h,\Z/p\Z)^* \\
&\cong&\ds\varinjlim_g \big(H^{i}(g/g\cap H,\Z/p\Z) \otimes H^{j}(h,\Z/p\Z)\big)^* \\
&\cong&\ds\varinjlim_g H^{i}\big(g/g\cap H, H^{j}(h,\Z/p\Z)\big)^*.
\end{array}
\]
Passing to the limit over $h$, we obtain the required isomorphism
\[
D_{i} (G/H,\Z/p\Z) \otimes D_{j}(H,\Z/p\Z) \cong D^2_{ij}(G,H,\Z/p\Z).
\]
\end{proof}

\begin{corollary}\label{lemma} Under the assumptions of Proposition~{\rm\ref{spectral}}, let $i_0$ and $ j_0$ be the smallest integers such that
$D_{i_0}(G/H,\Z/p\Z)\neq 0$ and $D_{j_0}(H,\Z/p\Z)\neq 0$, respectively. Then $D_{i_0+j_0}(G,\Z/p\Z)\neq 0$.
\end{corollary}

\begin{proof} The spectral sequence constructed in Proposition~\ref{spectral} induces an isomorphism
\[
D_{i_0+j_0}(G,\Z/p\Z) \cong D_{i_0}(G/H,\Z/p\Z) \otimes D_{j_0}(H,\Z/p\Z)\neq 0.
\]
\end{proof}

\begin{proof}[Proof of Theorem \ref{1}]
Assume that $G$ is a duality group at $p$ of dimension $d$. Let $\cd_p\, H=m$ and $n=d-m$.
Then there exists an open subgroup $H_1$ of $H$ such that $H^m(H_1,\Z/p\Z)\neq 0$.
Let $G_1$ be an open subgroup of $G$ such that $H_1=G_1\cap H$. Then $G_1$ is a duality group at $p$ of dimension $d$,
$\cd_p\, H_1=m$ and $G_1/H_1$ is an open subgroup of $G/H$. We consider the exact sequence
$$
1\lpfeil H_1\lpfeil G_1 \lpfeil G_1/H_1 \lpfeil 1.
$$
As $H^m(H_1,\Z/p\Z)$ is finite and nonzero, we have $\vcd_p\,G_1/H_1=n$, see \cite{NSW} (3.3.9). Furthermore,
$D_i(G_1,\Z/p\Z)=0$, $i<n+m$.
Using Corollary~\ref{lemma}, we see that $D_i(G_1/H_1,\Z/p\Z)=0$ for all $i<n$ and $D_j(H_1,\Z/p\Z)=0$  for all $j<m$.
Thus $G_1/H_1$, hence $G/H$, is virtually a duality group at $p$ of dimension $n$, and
 $H_1$, and so $H$, is a duality group at $p$ of dimension $m$. This shows (i).

Assume now that $H$ and $G/H$ are duality groups at $p$ of dimension $m$ and $n$.  Then, $\cd_p G=n+m$ by \cite{NSW} (3.3.8), and in the spectral sequence of Proposition~\ref{spectral} we have
 $E^2_{ij}=0$ for $(i,j)\neq (n,m)$. Hence $D_r(G,\Z/p\Z)=0$ for $r\neq n+m$ showing that $G$ is a duality group at $p$ of dimension $n+m$.

In  order to prove the assertion about the dualizing modules, let $h$ run through all open subgroups of $H$ which are normal in $G$ and $g$ runs through the open subgroups of $G$. Since $m=\cd_p\,H$, the Hochschild-Serre spectral sequence induces isomorphisms
\[
H^{m+n}(g,\Z/p^\nu\Z)\cong H^n(g/g\cap H,H^m(g\cap H,\Z/p^\nu\Z)),
\]
and we obtain
\renewcommand{\arraystretch}{1.5}

\[
\begin{array}{rcl}
I(G) & \cong& \ds\varinjlim_\nu\,\ds\varinjlim_g\,H^{m+n}(g,\Z/p\null^\nu\Z)^*  \\
     & \cong &\ds\varinjlim_\nu\,\ds\varinjlim_{h}\,
         \ds\varinjlim_{g}\,H^n\big(g/g\cap H,H^m(h,\Z/p\null^\nu\Z)\big)^*\\
     & \cong &\ds\varinjlim_\nu\,\ds\varinjlim_{h}\,
         \ds\varinjlim_{g,\res}\,H^0\big(g/g\cap H,\Hom\,(H^m(h,\Z/p\null^\nu\Z),I(G/H))\big) \\
     & \cong & \ds\varinjlim_\nu\,\ds\varinjlim_{h}\,\Hom(H^m(h,\Z/p\null^\nu\Z), I(G/H))\\
     &\cong &\Hom_{cts}(\ds\varprojlim_\nu\,\ds\varprojlim_{h}\,H^m(h,\Z/p\null^\nu\Z), I(G/H))\\
     &\cong &\Hom_{cts}\big((\ds\varinjlim_\nu\,\ds\varinjlim_{h} H^m(h,\Z/p^\nu\Z)^*)^\vee ,I(G/H)\big)\\
     & \cong & \Hom_{cts}\,(I(H)^\vee,I(G/H))\cong\big(I(H)^\vee\hatotimes_{\Z_p} I(G/H)^\vee\big)^\vee
\end{array}
\]
(see \cite {NSW} (5.2.9) for the last isomorphism).
This completes the proof of the theorem.
\end{proof}

\bigskip

\footnotesize{Alexander Schmidt, NWF I - Mathematik, Universit\"{a}t Regensburg, D-93040
Regensburg, Deutschland. email: alexander.schmidt@mathematik.uni-regensburg.de}

\bigskip
\footnotesize{Kay Wingberg, Mathematisches Institut, Universit\"{a}t Heidelberg, Im Neuenheimer Feld 288, 69120 Heidelberg, Deutschland. email: wingberg@mathi.uni-heidelberg.de}

\end{document}